\documentclass[12pt]{amsart}

\usepackage{a4wide}
\usepackage{latexsym}
\usepackage{amsfonts}
\usepackage{amsmath}
\usepackage{euscript}
\usepackage[all]{xy}
\usepackage{amscd}
\usepackage{amsthm}
\usepackage{amssymb}
\usepackage{verbatim}

\newtheorem{thm}{Theorem}[subsection]
\newtheorem{prop}[thm]{Proposition}
\newtheorem{defn}[thm]{Definition}

\newtheorem{lem}[thm]{Lemma}

\newtheorem{ex}[thm]{Example}

\def \Ob{\mathop{\rm Ob}}
\def \Hom{\mathop{\rm Hom}}
\def \Ext{\mathop{\rm Ext}}
\def \k{\underline{k}}
\def \End{\mathop{\rm End}}
\def \Aut{\mathop{\rm Aut}}
\def \Rad{\mathop{\rm Rad}}
\def\Mor{\mathop{\rm Mor}}

\def\Id{\mathop{\rm Id}}

\def\C{\mathcal{C}}

\def\A{\mathcal{A}}

\def\P{\mathcal{P}}

\def\lim{\mathop{\varprojlim}}
\def\colim{\mathop{\varinjlim}}
\def\H{\mathop{\rm H}}
\def\HH{\mathop{\rm HH}}

\def\Res{\mathop{\rm Res}}

\def\LG{\mathop{\Lambda\otimes_{k}\Gamma^{op}}}

\def\Z{\mathbb{Z}}
\def\E{\mathcal{E}}
\def\R{\mathcal{R}}
\begin{document}

\title{Hochschild and ordinary cohomology rings of small categories}

\author{Fei Xu}

\address{UMR 6629 CNRS/UN, Laboratoire de
Math\'ematiques Jean Leray, Universit\'e de Nantes, 2 Rue de la
Houssini\`ere, 44322 Nantes, France.}

\email{xu@math.univ-nantes.fr}

\maketitle

\begin{abstract} Let $\C$ be a small category and $k$ a field. There are two
interesting mathematical subjects: the category algebra $k\C$ and
the classifying space $|\C|=B\C$. We study the ring homomorphism
$\HH^*(k\C)\to\H^*(|\C|,k)$ and prove it is split surjective, using
the factorization category of Quillen \cite{Qu} and certain
techniques from functor cohomology theory. This generalizes the
well-known theorems for groups and posets. Based on this result, we
construct a seven-dimensional category algebra whose Hochschild
cohomology ring modulo nilpotents is not finitely generated,
disproving a conjecture of Snashall and Solberg \cite{SS}.\\

\end{abstract}

\pagenumbering{arabic}

\section{Introduction}

Let $\C$ be a small category, $k$ a field and $Vect_k$ the category
of $k$-vector spaces. We denote by $\Ob\C$ and $\Mor\C$ the sets of
objects and morphisms in $\C$, respectively. The category algebra
$k\C$ \cite{We, X1} of $\C$ is a $k$-vector space with basis
equal to $\Mor\C$, and the multiplication is given by the
composition of base elements (if two morphisms are not composable
then the product is zero). Suppose $Vect_k^{\C}$ is the category of
all covariant functors from $\C$ to $Vect_k$ and $k\C$-mod is the
category of left $k\C$-modules. Mitchell \cite[Theorem 7.1]{Mi}
showed that there exists a full faithful functor
$$
R: Vect_k^{\C} \to k\C\mbox{\rm{-mod}},
$$
defined by $R(F)=\oplus_{x\in\Ob\C} F(x)$. The functor $R$ has a
left inverse $L : k\C$-mod $\to Vect_k^{\C}$ defined by $M \mapsto
F_M$ such that $F_M(x)=1_x\cdot M$, where $1_x$ is the identity in
$\End_{\C}(x)$ for each $x\in\Ob\C$. When $\Ob\C$ is finite, the
category algebra $k\C$ has an identity
$1_{k\C}=\sum_{x\in\Ob\C}1_x$, and the above two functors provide an
equivalence between the two abelian categories. If $\C$ is a group
(regarded as a category with one object), the equivalence simply
gives us the fundamental correspondence between group modules and
group representations. In the present article we shall investigate
$\Ext^*_{Vect_k^{\C}}(M,N)=\oplus_{i\ge 0}\Ext^i_{Vect_k^{\C}}(M,N)$
for various $\C$ and functors $M,N \in Vect_k^{\C}$. Due to the
existence of the above faithful functor $R$, every functor is a
$k\C$-module. For simplicity, throughout this article we shall write
the above Ext as $\Ext^*_{k\C}(M,N)$. Whenever we need to emphasize that
a $k\C$-module $M$ is indeed an object in $Vect_k^{\C}$,
we say $M$ is a functor in $k\C$-mod. Let $\theta : \C_1 \to \C_2$
be a covariant functor between small categories. We use frequently
the functor $\Res_{\theta} : Vect_k^{\C_2} \to Vect_k^{\C_1}$,
which is called the restriction along $\theta$ (precomposition with $\theta$).
The functor $\theta$ does not always induce an algebra homomorphism
from $k\C_1$ to $k\C_2$ \cite{X1}. Hence it does not give rise to
a functor $k\C_2$-mod $\to k\C_1$-mod. Despite this potential
hole, in Section 2 we often write $\Res_{\theta} : k\C_2$-mod $\to k\C_1$-mod, again
for simplicity and consistency. As almost all modules we consider
are functors, it will not cause any real problem.

Let $\k \in k\C$-mod be the constant functor, sending every object
to $k$ and every morphism to the identity. When $\C$ is a group,
$\k=k$ becomes the trivial group module. For this reason, the
functor $\k$ is often called the trivial $k\C$-module, and it plays
the role of trivial module for a group algebra. The ordinary
cohomology ring of $\C$ with coefficients in $k$ can be defined as
$\Ext^*_{k\C}(\k,\k)$, which is isomorphic to $\H^*(|\C|,k)$
\cite{We, X1} and hence is graded commutative. Such an
ordinary cohomology ring modulo nilpotents is not finitely generated
in general, see for example \cite{X2}.

Let $\C^e=\C\times\C^{op}$, where $\C^{op}$ is the opposite
category. The enveloping algebra of $k\C$,
$(k\C)^e=k\C\otimes_k(k\C)^{op}$, is naturally isomorphic to $k\C^e$
as $k$-algebras. Hence in the present article we shall not
distinguish the two algebras $(k\C)^e$ and $k\C^e$. By introducing
$\C^e$ and $k\C^e$, one can use functor cohomology theory to
investigate Hochschild cohomology. We want to consider
${\Ext}^*_{k\C^e}(M,N)$, where $M, N\in k\C^e$-mod. When $M=N=k\C$,
$\Ext^*_{k\C^e}(k\C,k\C)$ becomes a graded commutative ring \cite{SS}. If $\Ob\C$ is finite (thus $k\C$
has an identity), one can identify the above ring with the Hochschild cohomology ring $\HH^*(k\C)$ (see
\cite[Section 7]{Ge} and \cite[Chapter 1]{Lo}). For this reason, we
shall call $\Ext^*_{k\C^e}(k\C,k\C)$ the Hochschild
cohomology ring of $\C$ in the present article. We note that the
module $k\C\in k\C^e$-mod comes from a functor $\C^e \to Vect_k$
such that $k\C(x,y)=k\Hom_{\C}(y,x)$ for each $(x,y)\in\Ob\C^e$
(if $\Hom_{\C}(y,x)=\emptyset$ then we assume $k\C(x,y)=0$).

Suppose $A$ is an associative $k$-algebra and $A^e$ is its
enveloping algebra. Let $M$ be an $A^e$-module. Then one has a ring
homomorphism induced by the tensor product $-\otimes_A M$
$$
\phi_M : {\Ext}^*_{A^e}(A,A) \to {\Ext}^*_A(M,M).
$$
If we take $A=k\C$ for a small category $\C$ and $M=\k$, we get a ring homomorphism
$$
{\phi}_{\C} : {\Ext}^*_{k\C^e}(k\C,k\C) \to {\Ext}^*_{k\C}(\k,\k).
$$
In this situation, ${\phi}_{\C}$ is really induced by the projection
functor $pr : \C^e \to \C$ (see Section 2.3). The structures of these
two cohomology rings and the homomorphism are the main subjects of
our investigation. Note that we name the ring homomorphism
$\phi_{\C}$, not $\phi_{\k}$, since we need to deal with various
categories and $\phi_{\k}$ can cause confusion. It is well-known
that when $\C$ is a group, $\phi_{\C}$ is a split surjection (see
for instance \cite{Be} or \cite[2.9]{L1}), whilst, when $\C$ is a poset,
$\phi_{\C}$ is an isomorphism \cite{GS}. The two results are proved
in completely different ways in the literature. In our article, we
use functor cohomology theory to establish a general statement on
the ring homomorphism $\phi_{\C}$, including the above two results
as special cases. In order to deal with the general situation, we
need to consider the category of factorizations in a category $\C$,
introduced by Quillen \cite{Qu}. The category of factorizations in
$\C$, $F(\C)$, has all the morphisms in $\C$ as its objects. If we
write the objects in $F(\C)$ as $[\alpha]$, for any
$\alpha\in\Mor\C$, then there exists morphisms from $[\alpha]$ to
$[\alpha']$ if $\alpha$ factors through $\alpha'$ in $\Mor\C$. The
category $F(\C)$ admits natural functors $t$ and $s$ into $\C$ and
$\C^{op}$, respectively, inducing homotopy equivalences of classifying spaces. One can
assemble these two functors together to form a new functor
$\tau=(t,s) : F(\C) \to \C^e$. Quillen observed that $F(\C)$ is
cofibred over $\C^e$ and described the fibres. Based on these, we
prove the following statements
(Theorem 2.3.1 and Proposition 2.3.5). We comment that Mac Lane \cite{ML}
discussed the question for monoids in Section X.5 of his book and obtained
part of the result (stated for homology).\\

\noindent{\bf Theorem A} {\it Let $\C$ be a small category and $k$ a
field. For any functor $M \in k\C^e$-mod, we have
$$
{\Ext}^*_{k\C^e}(k\C,M)\cong{\Ext}^*_{kF(\C)}(\k,{\Res}_{\tau}M),
$$
where $\Res_{\tau}$ is the restriction along $\tau : F(\C)\to\C^e$
(precomposition with $\tau$). In particular we have
$$
{\Ext}^*_{k\C^e}(k\C,\k)\cong{\Ext}^*_{kF(\C)}(\k,\k)\cong{\Ext}^*_{k\C}(\k,\k),
$$
and $\phi_{\C} : {\Ext}^*_{k\C^e}(k\C,k\C)\to{\Ext}^*_{k\C}(\k,\k)$ is a split surjection, induced by the following
decompositions $\Res_{\tau}(k\C)\cong\k\oplus N_{\C}$ and
$$
{\Ext}^*_{k\C^e}(k\C,k\C)\cong{\Ext}^*_{k\C}(\k,\k)\oplus{\Ext}^*_{kF(\C)}(\k,N_{\C}),
$$
where $N_{\C} \in kF(\C)$-mod as a functor takes the following value
$$
N_{\C}([\alpha])=k\{\beta - \gamma \bigm{|} \beta, \gamma \in
{\Hom}_{\C}(y,x) \},
$$
if $[\alpha]\in\Ob F(\C)$ and $\alpha\in\Hom_{\C}(y,x)$.}\\

Especially, the existence of a surjective homomorphism implies that
if the ordinary cohomology ring, modulo nilpotents, is not finitely
generated, neither is the Hochschild cohomology ring, modulo
nilpotents. In \cite{X2} we computed the mod-2 ordinary cohomology
ring of the following category $\E_0$
$$
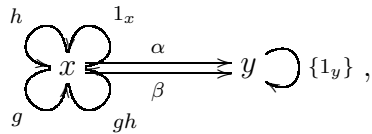
\xymatrix{x \ar@(r,u)[]_{1_x} \ar@(l,d)[]_g \ar@(u,l)[]_h
\ar@(d,r)[]_{gh} \ar@<2pt>[rr]^{\alpha} \ar@<-2pt>[rr]_{\beta}&& y
\ar@(ur,dr)[]^{\{1_y\}}} ,
$$
where $g^2=h^2=1_x, gh=hg, \alpha h = \beta g = \alpha$, and $\alpha
g = \beta h = \beta$. It was shown there that its ordinary
cohomology ring doesn't have any nilpotents and is not finitely
generated. Thus its Hochschild cohomology ring modulo nilpotents is
not finitely generated, providing a counterexample against the
conjecture in \cite{SS}. We note that the category algebra $k\C$ is
not a self-injective algebra, in contrast to the fact
that the Hochschild cohomology ring of a finite-dimensional
cocommutative Hopf algebra, or of a finite-dimensional
self-injective algebra of finite representation type, is finitely
generated \cite{FS, GSS} (a finite-dimensional Hopf algebra is always
self-injective \cite{LS}). In this particular case, the
category algebra is graded and is Koszul, which was brought attention to the
author by Nicole Snashall.

A small category is called EI if every endomorphism is an
isomorphism. A category if finite if the morphism set is finite.
Typical examples of finite EI-categories are posets and groups. The
above category $\E_0$ is finite EI as well. Some sophisticated
finite EI-categories have been heavily used in, for example, the
$p$-local finite group theory \cite{BLO} and modular representation
theory \cite[Chapter 7]{Th}. Let $\C$ be a finite EI-category. We
can define a full subcategory $\A_{\C}=\A$ such that $\Ob\A=\Ob\C$
and $\Mor\A$ contains exactly all the isomorphisms in $\Mor\C$. The
category $\A$ can be considered as the disjoint union of all
finite groups in $\C$. The following is Theorem 2.4.2.\\

\noindent{\bf Theorem B} Let $\C$ be a finite EI-category and $k$ a
field. Then we have the following commutative diagram
$$
\xymatrix{{\Ext}^*_{k\C^e}(k\C,k\C) \ar[r]^{\phi_{k\A}} \ar[d]_{\phi_{\C}} & {\Ext}^*_{k\A^e}(k\A,k\A) \ar[d]^{\phi_{\A}}\\
{\Ext}^*_{k\C}(\k,\k) \ar[r]_{\Res_{\C,\A}} &
{\Ext}^*_{k\A}(\k,\k).}
$$
Here $\Res_{\C,\A}$ is induced by the inclusion $\iota : \A
\hookrightarrow \C$. In this theorem the category $\A$ may be
replaced by any full subcategory of it.

Our paper begins with a brief introduction to the ring homomorphisms
from the Hochschild cohomology of an associative algebra to some
relevant rings. Afterwards, we introduce the concept of an
enveloping category and reinterpret the ring homomorphism using
functor cohomology theoretic methods. Based on Quillen's work, we
continue to prove $\phi_{\C} :
{\Ext}^*_{k\C^e}(k\C,k\C)\to{\Ext}^*_{k\C}(\k,\k)$ is split
surjective for any small category $\C$. Some consequences of this
splitting surjection and further properties will be given. Finally,
we end this paper with four examples. The first example provides a
counter-example to a conjecture of Snashall
and Solberg.\\

\noindent{\bf Acknowledgements} I wish to thank Aur\'elien Djament,
Laurent Piriou and Vincent Franjou, my co-investigators working on
the CNRS research project ``Functor Homology Theory'', for many
helpful conversations. The counterexample against the finite
generation of ordinary cohomology rings of finite EI-categories (see
$\E_0$ in Section 3.1), which we found together, was one of the
starting points of the present paper. I also would like to thank
Nicole Snashall for useful comments on my work. The author is
partially supported by a CNRS post-doctoral fellowship.

\section{Hochschild and ordinary cohomology rings of categories}

We first describe the ring homomorphism from the Hochschild
cohomology ring, of an associative algebra, to some relevant
cohomology rings, induced by tensor products with modules. When the
associative algebra is a category algebra and the target is the ordinary
cohomology ring, we reconstruct the ring
homomorphism, using a different method. Based on the alternative
description, we show the ring homomorphism $\phi_{\C}$ is split surjective.

\subsection{The ring homomorphisms from the Hochschild cohomology ring}

\begin{defn} Let $A$ be an associative $k$-algebra and $M, N$ two $A$-modules.
We write $\Ext^*_A(M,N) = \oplus_{i\ge 0} \Ext^i_A(M,N)$.
\end{defn}

In general, if $\Lambda$ and $\Gamma$ are two
associative $k$-algebras and $M$ is a
$\Lambda\otimes_k\Gamma^{op}$-module, or equivalently a
$\Lambda$-$\Gamma$-bimodule, we can define a ring homomorphism
induced by the tensor product $-\otimes_{\Lambda}M$
$$
\phi_M : {\Ext}^*_{\Lambda^e}(\Lambda, \Lambda) \to
{\Ext}^*_{\Lambda\otimes_k\Gamma^{op}}(M, M).
$$
Let $\R_* \to \Lambda \to 0$ be a projective resolution of the
$\Lambda^e$-module $\Lambda$. The exact sequence is split if we
regard it as a complex of right $\Lambda$-modules. Thus by tensoring
$M$ over $\Lambda$ from the right, we obtain an exact sequence
ending at the $\Lambda\otimes_k\Gamma^{op}$-module $M$
$$
\R_*\otimes_{\Lambda}M \to \Lambda\otimes_{\Lambda}M \cong M \to 0.
$$
Now one can build a projective resolution of $M$, $\R'_*\to M \to
0$, along with a chain map
$$
\xymatrix{\R'_* \ar[r] \ar[d] & M \ar[r] \ar[d]^{=} & 0\\
\R_*\otimes_{\Lambda}M \ar[r] & \Lambda\otimes_{\Lambda}M \ar[r] &
0.}
$$
This induces an algebra homomorphism $\phi_M :
\Ext^*_{\Lambda^e}(\Lambda,\Lambda) \to
\Ext^*_{\Lambda\otimes_k\Gamma^{op}}(M,M)$. If $N$ is another
$\Lambda\otimes_k\Gamma^{op}$-module, we see
$\Ext^*_{\Lambda\otimes_k\Gamma^{op}}(M,N)$ has an
$\Ext^*_{\Lambda^e}(\Lambda,\Lambda)$-module structure via the ring
homomorphisms $\phi_M$ and $\phi_N$ together with the Yoneda splice.
We quote the following theorem of Snashall and Solberg \cite{SS}.

\begin{thm} Let $\Lambda$ and $\Gamma$ be two associative $k$-algebras. Let $\eta$
be an element in $\Ext^n_{\Lambda^e}(\Lambda,\Lambda)$ and $\theta$
an element in $\Ext^m_{\LG}(M,N)$ for two
$\Lambda$-$\Gamma$-bimodules $M$ and $N$. Then
$\phi_N(\eta)\theta=(-1)^{mn}\theta\phi_M(\eta)$.
\end{thm}

When $\Lambda$ has an identity, it means
$\Ext^*_{\Lambda^e}(\Lambda,\Lambda)\cong\HH^*(\Lambda)$ is a graded
commutative ring, which was first proved by Gerstenhaber \cite{Ge}.

\subsection{Enveloping category of a small category}

Let $\C$ be a small category. Quillen \cite[page 94 Example]{Qu}
considered the category $\C^{op}\times\C$. We slightly modify it and
give it a name, in order to be consistent with our investigation of
the Hochschild cohomology.

\begin{defn}
We call $\C^e = \C \times \C^{op}$ the enveloping category of a
small category $\C$.
\end{defn}

The following result is just a simple observation. It implies the
enveloping algebra of a category algebra of $\C$ is the category
algebra of its enveloping category, so later on we will just use the
terminology $k\C^e$ when dealing with Hochschild cohomology. This
identification enables us to apply functor cohomology theory to the
investigation of the Hochschild cohomology theory of category
algebras.

\begin{lem} Let $\C$ be a small category. There is a natural isomorphism
$k\C^e\cong(k\C)^e$. As a functor, $k\C(x,y)=k\Hom_{\C}(y,x)$ if $\Hom_{\C}(y,x)\ne\emptyset$
and $k\C(x,y)=0$ otherwise. Here $(x, y)\in\Ob\C^e$.

\begin{proof}[Proof.] We define a map $k\C^e \to (k\C)^e$ on the natural base elements
of $k\C^e$ by $(\alpha,\beta^{op})\mapsto\alpha\otimes\beta^{op}$,
$\alpha, \beta \in \Mor\C$. It extends linearly to an algebra
isomorphism.

If $M$ is a $k\C^e$-module and $m\in M$, then
$(\alpha,\beta^{op})\cdot m = \alpha\cdot m \cdot \beta$ and as a
functor $M : \C^e \to Vect_k$
$$
M(x,y)=1_{(x,y)}\cdot M=(1_x,1_y^{op})\cdot M = 1_x\cdot M \cdot 1_y,
$$
on each object $(x,y)\in\C^e$. In particular,
$$
k\C(x,y)=(1_x,1_y^{op})\cdot k\C=1_x\cdot k\C \cdot 1_y =
k{\Hom}_{\C}(y,x)
$$
if $\Hom_{\C}(y,x)\ne\emptyset$, and $k\C(x,y)=0$ otherwise.
\end{proof}
\end{lem}

Let $\C$ be a small category. We recall Quillen's category $F(\C)$
of factorizations in $\C$. In his article \cite{Qu}, Quillen named this category $S(\C)$.
However since $S(\C)$ has been used to denote the subdivision of a small category
$\C$ \cite{Sl, L2}, we adopt Baues and Wirsching's terminology \cite{BW}
which we believe is suitable. The category $F(\C)$ has the morphisms in
$\C$ as its objects. In order to avoid confusion, we write an object
in $F(\C)$ as $[\alpha]$, whenever $\alpha\in\Mor\C$. A morphism from
$[\alpha]\in\Ob F(\C)$ to $[\alpha']\in\Ob F(\C)$ is given by a pair
of $u, v \in \Mor\C$, making the following diagram commutative
$$
\xymatrix{x \ar[d]_u & y \ar[l]_{\alpha} \ar[d]^{v^{op}}\\
x' & y'. \ar[l]^{\alpha'}}
$$
In other words, there is an morphism from $[\alpha]$ to $[\alpha']$
if and only if $\alpha'=u\alpha v$ for some $u,v\in\Mor\C$, or
equivalently $\alpha$ is a factor of $\alpha'$ in $\Mor\C$. The category $F(\C)$
admits two natural covariant functors to $\C$ and $\C^{op}$
$$
\xymatrix{\C & F(\C) \ar[l]_t \ar[r]^s & \C^{op}},
$$
where $t$ and $s$ send an object $[\alpha]$ to its target and
source, respectively. Using his Theorem A and its corollary, Quillen
showed these two functors induce homotopy equivalences of the
classifying spaces. We will be interested in the functor
$$
\tau=(t,s) : F(\C) \to \C^e = \C\times\C^{op},
$$
sending an $[\alpha] \in \Ob F(\C)$ to $(x,y) \in \Ob\C^e$
if $\alpha\in\Hom_{\C}(y,x)$ and a morphism $(u,v^{op})\in\Mor
F(\C)$ to $(u,v^{op})\in\Mor(\C^e)$.

The importance of the functor $\tau : F(\C) \to \C^e$ lies in the
fact that its target category gives rise to the Hochschild
cohomology ring of $\C$, while its source category determines the
ordinary cohomology ring of $\C\simeq F(\C)$. In the situation of
(finite) posets and groups, the functor is well-understood and in
the group case it has been implicitly used to establish the
homomorphism from the Hochschild cohomology ring to the ordinary
cohomology ring.

\begin{ex}
\begin{enumerate}
\item When $\C$ is a poset, $\tau : F(\C) \to \C^e$ sends $F(\C)$ isomorphically
onto a full category $\C^e_{\Delta}\subset\C^e$, where
$$
{\Ob\C}^e_{\Delta}=\{(x,y)\in\Ob\C^e\bigm{|}{\Hom}_{\C}(y,x)\ne\emptyset\}
$$
(the full subcategory $\C^e_{\Delta}$ is well-defined whenever $\C$
is EI, see Section 2.4). One can easily see that $k\C$ as a functor
only takes non-zero values at objects in $\Ob\C^e_{\Delta}$.
Furthermore as a $k\C^e_{\Delta}$-module, $k\C \cong \k$ is the
trivial module by Lemma 2.2.2. Since $\C^e_{\Delta} \cong F(\C)$ is
a co-ideal in the poset $\C^e$, we obtain
$\Ext^*_{k\C^e}(k\C,k\C)\cong\Ext^*_{k\C^e_{\Delta}}(k\C,k\C)\cong\Ext^*_{kF(\C)}(\k,\k)\cong\Ext^*_{k\C}(\k,\k)$,
where the last isomorphism comes from the fact that $|F(\C)| \simeq
|\C|$. This isomorphism between the two cohomology rings was first
established in \cite{GS};

\item When $\C$ is a group, the category
$F(\C)$ is a groupoid and is equivalent to a subcategory of the one
object category $\C^e$ with morphism set
$$
\{(g,{g^{-1}}^{op})\bigm{|} g\in\Mor\C\}\subset\Mor\C^e.
$$
Based on this description, one can prove the existence of the surjective homomorphism from the
Hochschild cohomology ring to the ordinary cohomology ring of a
group, which is basically the same as the classical approach.
See for example \cite{Be}.
\end{enumerate}
\end{ex}

\subsection{The main theorem} In order to deal with the general situation, we need to recall the
definition of an overcategory. It is used to define and understand
the left Kan extension, which generalizes the concept of an
induction.

Let $\theta : \C_1 \to \C_2$ be a covariant functor between small
categories. For each $z\in\Ob\C_2$, the overcategory $\theta/z$
consists of objects $(x,\alpha)$, where $x\in\Ob\C_1$ and
$\alpha\in\Hom_{\C_2}(\theta(x),z)$. A morphism from $(x,\alpha)$ to
$(x',\alpha')$ is a morphism $\beta \in \Hom_{\C_1}(x,x')$ such that
$\alpha=\alpha'\theta(\beta)$. Let $\Res_{\theta} : k\C_2$-mod $\to
k\C_1$-mod be the restriction on functors along $\theta$ (precomposition with
$\theta$). The left adjoint of $\Res_{\theta}$ is called the left
Kan extension $LK_{\theta} : k\C_1$-mod $\to k\C_2$-mod and is
defined by
$$
{LK}_{\theta}(M)(z)={\colim}_{\theta/z}M\circ\pi,
$$
where $z\in\Ob\C_2$, $\pi : \theta/z \to \C_1$ is the projection
functor $(x,\alpha) \mapsto x$ and $M$ is a functor in $k\C_1$-mod. When
$\C_2$ is a subgroup of a group $\C_1$ and $\theta$ is the
inclusion, the left Kan extension is the usual induction, i.e.
$LK_{\theta}(M)\cong k\C_2\otimes_{k\C_1}M$.

With the definition of an overcategory, one can continue to define
two functors $\theta/? : \C_2 \to$ sCat (the category of small
categories), and $C_*(\theta/?) : \C_2 \to k\C_2$-Cplx (the category
of complexes of $k\C_2$-modules). For each $x\in\Ob\C_2$,
$C_*(\theta/x)$ is the simplicial complex coming from the nerve of
the small category $\theta/x$. When $\C_1=\C_2=\C$ and
$\theta=\Id_{\C}$, we have functors $\Id_{\C}/?$ and
$C_*(\Id_{\C}/?)$. It is well-known that the latter can be used to
define a projective resolution of the $k\C$-module $\k$ :
$C_*(\Id_{\C}/?) \to \k \to 0$. For each $n \ge 0$, $C_n(\Id_{\C}/?)
: \C \to k\C$-Cplx is the functor sending each $x\in\Ob\C$ to the
vector space whose basis is the set of all $n$-chains of morphisms
in $\Id_{\C}/x$. The differential, a $k\C$-map, $\sigma^n :
C_n(\Id_{\C}/?) \to C_{n-1}(\Id_{\C}/?)$ is defined as follows. For
each $x \in \Ob\C$,
$$
\begin{array}{ll}
&\sigma^n_x((x_0,\alpha_0) \to \cdots \to (x_i,\alpha_i) \to \cdots
\to
(x_n,\alpha_n))\\
=& \sum^n_{i=0}(-1)^i[(x_0,\alpha_0) \to \cdots \to
\widehat{(x_i,\alpha_i)} \to \cdots \to (x_n,\alpha_n)],
\end{array}
$$
where $\alpha_i \in \Hom_{\C}(x_i,x)$. Let $\theta : \C_1 \to \C_2$
be a covariant functor. There is an isomorphism of complexes of
projective $k\C_2$-modules (a left Kan extension always preserves
projectives)
$$
LK_{\theta}(C_*({\Id}_{\C_1}/?)) \cong C_*(\theta/?),
$$
which can be found for example in Hollender-Vogt \cite[4.3]{HV}. Under certain
conditions, the above complex may be a
projective resolution of the $k\C_2$-module $LK_{\theta}(\k)$. This
is the key to our future investigation.

We want to discuss the left Kan extensions of the functors $\tau$,
$t$ and $pr$ in the following commutative diagram of small categories
$$
\xymatrix{F(\C) \ar[dr]_t \ar[rr]^{\tau} & & \C^e=\C\times\C^{op} \ar[dl]^{pr}\\
&\C&,}
$$
where $pr$ is the projection onto the first component. Since $t=pr\circ\tau$, we have
$$
{LK}_t \cong {LK}_{pr}\circ{LK}_{\tau}.
$$
In the rest of this section, we will establish and describe the following ring
homomorphisms, induced by the three left Kan extensions $LK_t,
LK_{pr}$ and $LK_{\tau}$ respectively,
$$
\begin{array}{lllll}
t^* & : & {\Ext}^*_{kF(\C)}(\k,\k) & \to & {\Ext}^*_{k\C}(\k,\k),\\
pr^* & : & {\Ext}^*_{k\C^e}(k\C,k\C) & \to & {\Ext}^*_{k\C}(\k,\k)\\
\tau^* & : & {\Ext}^*_{kF(\C)}(\k,\k) & \to &
{\Ext}^*_{k\C^e}(k\C,k\C).
\end{array}
$$
The first two homomorphisms are not difficult to describe and we do
it now. The homomorphism $t^*$ is an isomorphism since $t$ induces a
homotopy equivalence of $F(\C)$ and $\C$ by \cite{Qu}. More
explicitly, let $C_*(\Id_{F(\C)}/?) \to \k \to 0$ be the projective
resolution of the $kF(\C)$-module $\k$. The left Kan extension of
$t$, $LK_t$, sends it to a projective resolution of the $k\C$-module
$\k$
$$
{LK}_t(C_*({\Id}_{F(\C)}/?))\cong C_*(t/?) \to {LK}_t(\k)\cong\k \to 0.
$$
The reason is that first of all, $C_*(t/?)$ is a complex of
projective $k\C$-modules, and second of all, for each $x\in\Ob\C$,
$t/x$ is contractible \cite{Qu} and thus $C_*(t/x)$ is exact except
having homology $k$ at the end.

The homomorphism $pr^*$, induced by $pr$, is exactly $\phi_{\C}$, defined earlier,
which is induced by tensoring over $k\C$ with $\k$ from the right.
We see this from the fact that $LK_{pr}$ is exactly the tensor
product $-\otimes_{k\C}\k$ on a projective resolution of the
$k\C^e$-module $k\C$. In fact for each $x\in\Ob\C$ since $pr/x \cong
(\Id_{\C}/x) \times \C^{op}$,
$$
{LK}_{pr}(k\C^e)(x)={\colim}_{pr/x}k\C^e\cong{\colim}_{\Id_{\C}/x}(k\C)\otimes_k{\colim}_{\C^{op}}(k\C^{op})\cong
1_x\cdot k\C\otimes_k \k.
$$
It implies $LK_{pr}(k\C^e)\cong k\C\otimes_k \k\cong
k\C^e\otimes_{k\C}\k$. Also we have
$$
{LK}_{pr}(k\C)\cong{LK}_{pr}({LK}_{\tau}(\k))\cong{LK}_{t}(\k)\cong\k.
$$

Now we turn to investigate $LK_{\tau}$ and $\tau^*$. Our goal is to use
$\tau^*$ and $t^*$ to interpret $pr^*=\phi_{\C}$. The main result in
this section is as follows.

\begin{thm} Let $\C$ be a small category and $k$ a field. There exists a ring homomorphism
$$
\epsilon^* : {\Ext}^*_{k\C^e}(k\C,k\C) \to {\Ext}^*_{kF(\C)}(\k,\k)
$$
such that $\epsilon^*\tau^*\cong 1$. Moreover the following
composition $t^*\epsilon^*$ is a split surjection
$$
{\Ext}^*_{k\C^e}(k\C,k\C){\buildrel{\epsilon^*}\over{\twoheadrightarrow}}{\Ext}^*_{kF(\C)}(\k,\k){\buildrel{t^*}\over{\to}}{\Ext}^*_{k\C}(\k,\k),
$$
with the property that $t^*\epsilon^*\cong pr^*\cong\phi_{\C}$.
\end{thm}

The proof of this theorem will be divided into three lemmas. We
first discuss the action of $LK_{\tau}$ on a certain projective
resolution of the $kF(\C)$-module $\k$. In his example on page 94 of
\cite{Qu}, Quillen asserted that the category $F(\C)$ is a cofibred
category over $\C^e$, via $\tau$, with discrete fibres defined by
the functor $(x,y) \mapsto \Hom_{\C}(y,x)$, where $(x,y) \in
\Ob\C^e$. As a consequence of the assertion Quillen indicated that
each overcategory $\tau/(x,y)$ is homotopy equivalent to the fibre
$\tau^{-1}(x,y)$, which is the discrete category $\Hom_{\C}(y,x)$. Hence the
left Kan extension of $\k$ takes the following value at each object
$(x,y)$
$$
LK_{\tau}(\k)(x,y)={\colim}_{\tau/(x,y)}\k\cong{\H}_0(|\tau/(x,y)|,k)\cong{\H}_0(|\tau^{-1}(x,y)|,k),
$$
which equals $k{\Hom}_{\C}(y,x)$ if $\Hom_{\C}(y,x)\ne\emptyset$ and zero otherwise. It implies $LK_{\tau}(\k)\cong k\C$ as $k\C^e$-modules. Further
more, the following lemma implies $LK_{\tau}(C_*(\Id_{F(\C)}/?))\to
LK_{\tau}(\k)\cong k\C \to 0$ is indeed a projective resolution.

\begin{lem} Let $\C_1$ and $\C_2$ be two small categories and $\theta : \C_1 \to \C_2$
a covariant functor. If $\theta/w$ is a discrete
category for every $w\in\Ob\C_2$, then we obtain a projective
resolution of the $k\C_2$-module
$LK_{\theta}(\k)\cong\H_0(|\theta/?|,k)$
$$
LK_{\theta}(C_*({\Id}_{\C_1}/?)) \cong C_*(\theta/?) \to
LK_{\theta}(\k)\to 0.
$$

\begin{proof}[Proof.] Evaluating $C_*(\theta/?)$ at an object $w\in\Ob\C_2$, one gets a complex
$C_*(\theta/w)$ that computes the homology of $|\theta/w|$ with coefficients in $k$. Thus if $\theta/w$ is a discrete category, we get an
exact sequence
$$
LK_{\theta}(C_*({\Id}_{\C_1}/w)) \cong C_*(\theta/w) \to
LK_{\theta}(\k)(w)\cong{\H}_0(|\theta/w|,k) \to 0.
$$
If $\theta/w$ is a discrete category for every $w\in\Ob\C_2$, then
we obtain a projective resolution of the $k\C_2$-module
$LK_{\theta}(\k)$
$$
LK_{\theta}(C_*({\Id}_{\C_1}/?)) \cong C_*(\theta/?) \to
LK_{\theta}(\k) \to 0,
$$
because it's exact and meanwhile the left Kan extension preserves
projectives.
\end{proof}
\end{lem}

Since $LK_{\theta}$ is the left adjoint of $\Res_{\theta}$, there
are natural transformations $\Id \to \Res_{\theta}LK_{\theta}$ and
$LK_{\theta}\Res_{\theta}\to\Id$. We pay attention to the case of
$\tau : F(\C)\to\C^e$. There exists a $kF(\C)$-homomorphism
$\k\to\Res_{\tau}LK_{\tau}(\k)=\Res_{\tau}(k\C)$ as well as a
$k\C^e$-homomorphism $k\C=LK_{\tau}\Res_{\tau}(\k)\to\k$. The latter
gives rise to a $kF(\C)$-homomorphism
$\Res_{\tau}(k\C)=\Res_{\tau}LK_{\tau}\Res_{\tau}(\k)\to\k=\Res_{\tau}\k$.
In case $\C$ is a poset, one has $\k=\Res_{\tau}(k\C)$. When
$\C$ is a group, $F(\C)$ is a groupoid, equivalent to the
automorphism group of $[1_{\C}]\in\Ob F(\C)$, that is,
$\{(g,{g^{-1}}^{op})\bigm{|} g\in\Mor\C\}$. If we name the full
subcategory of $F(\C)$, consisting of one object $[1_{\C}]$, by
$\tilde\Delta\C$ and the inclusion (an equivalence) by $i :
\tilde\Delta\C \hookrightarrow F(\C)$. Then $\Res_{\tau
i}(k\C)=\Res_{\tau}(k\C)([1_{\C}])$ is a $k\tilde\Delta\C$-module
with the action $(g,{g^{-1}}^{op})\cdot a = gag^{-1}$,
$a\in\Res_{\tau i}(k\C)$. Thus $\Res_{\tau i}(k\C)=\oplus kc_g$,
where $c_g$ is the conjugacy class of $g\in\Mor\C$. In particular
$\k = kc_{1_{\C}}$ is a direct summand of $\Res_{\tau i}(k\C)$ and
it implies $\k \bigm{|} \Res_{\tau}(k\C)$ as $kF(\C)$-modules
because $i$ is an equivalence of categories.

\begin{lem} Let $\C$ be a small category. Then
$\k \bigm{|} \Res_{\tau}(k\C)$ as $kF(\C)$-modules.

\begin{proof}[Proof.] One needs to keep in mind that the restriction
of a module usually has a large $k$-dimension than the module itself
since $\tau$ is not injective on objects. We define a
$kF(\C)$-homomorphism (a natural transformation) $\iota : \k \to
\Res_{\tau}(k\C)$ by the assignments  $\iota_{[\alpha]}
(1_k)=\alpha\in\Res_{\tau}(k\C)([\alpha])$ for each $[\alpha]\in\Ob
F(\C)$. If $[\beta]$ is another object in $\Ob F(\C)$ and
$(u,v^{op}) \in \Hom_{F(\C)}([\alpha],[\beta])$ is an arbitrary
morphism, then by the definition of an $F(\C)$-morphism,
$(u,v^{op})\cdot\alpha=u\alpha v =\beta$. Hence $\iota$ maps $\k$
isomorphically onto a submodule of $\Res_{\tau}(k\C)$. On the other
hand, we may define a $kF(\C)$-homomorphism $\epsilon :
\Res_{\tau}(k\C) \to \k$ such that, for any $[\alpha]\in\Ob F(\C)$,
$\epsilon_{[\alpha]} : \Res_{\tau}(k\C)([\alpha]) \to
\k([\alpha])=k$ sends each base element in
$\Res_{\tau}(k\C)([\alpha])=k\Hom_{\C}(y,x)$ to $1_k$. One can
readily check the composite of these two maps is the identity
$$
\k {\buildrel{\iota}\over{\to}} {\Res}_{\tau}(k\C)
{\buildrel{\epsilon}\over{\to}} \k,
$$
and this means $\k \bigm{|} \Res_{\tau}(k\C)$ or $\Res_{\tau}(k\C) =
\k\oplus N_{\C}$ for some $kF(\C)$-module $N_{\C}$.
\end{proof}
\end{lem}

The module $N_{\C}$ as a functor can be described by
$$
N_{\C}([\alpha])=k\{\beta - \gamma \bigm{|} \beta, \gamma \in
{\Hom}_{\C}(y,x) \},
$$
if $[\alpha]\in\Ob F(\C)$ and $\alpha\in\Hom_{\C}(y,x)$. It will be
useful to our computation since it determines the ``difference''
between the Hochschild and ordinary cohomology rings of a category.
The next lemma finishes off our proof of the main theorem.

\begin{lem} Let $\C$ be a small category. There is a surjective ring homomorphism $\epsilon^*$
$$
{\Ext}^*_{k\C^e}(k\C,k\C)\twoheadrightarrow{\Ext}^*_{kF(\C)}(\k,\k),
$$
such that $\epsilon^*\tau^*\cong 1$ and $pr^*\cong t^*\epsilon^*$.

\begin{proof}[Proof.] By Quillen's observation \cite{Qu}, we know every overcategory
$\tau/(x,y)$ has the homotopy type of $\Hom_{\C}(y,x)$. Applying
Lemma 2.3.2 to $\tau : F(\C) \to \C^e$, we know the left
Kan extension $LK_{\tau}$ sends a certain projective resolution
$\P_*$ of the $kF(\C)$-module $\k$ to a projective resolution
$LK_{\tau}(\P_*)$ of the $k\C^e$-module $k\C$. Then on the cochain level we see $\tau^*$ is
determined by the following composition.
$$
{\Hom}_{kF(\C)}(\P_*,\k) \to
{\Hom}_{k\C^e}(LK_{\tau}(\P_*),LK_{\tau}(\k)) \cong
{\Hom}_{kF(\C)}(\P_*,{\Res}_{\tau}LK_{\tau}(\k)).
$$
Lemma 2.3.3 says ${\Res}_{\tau}LK_{\tau}(\k)=\k\oplus N_{\C}$ for
some $kF(\C)$-module $N_{\C}$. As a consequence, we have a split
exact sequence of $k$-vector spaces
$$
0\to{\Ext}^*_{kF(\C)}(\k,\k) \hookrightarrow
{\Ext}^*_{k\C^e}(k\C,k\C) \cong {\Ext}^*_{kF(\C)}(\k,\k\oplus
N_{\C}) \twoheadrightarrow{\Ext}^*_{kF(\C)}(\k,\k)\to 0.
$$
The leftmost map is $\tau^*$ and the rightmost map is named
$\epsilon^*$, induced by $\epsilon$ in Lemma 2.3.3, and is given by
$$
{\Hom}_{k\C^e}(LK_{\tau}(\P_*),LK_{\tau}(\k)) \cong
{\Hom}_{kF(\C)}(\P_*,{\Res}_{\tau}LK_{\tau}(\k))\to{\Hom}_{kF(\C)}(\P_*,\k).
$$
From here, we can see $pr^*\cong t^*\epsilon^*$ because of the
following commutative diagram
$$
\xymatrix{\Hom_{k\C^e}(LK_{\tau}(\P_*),LK_{\tau}(\k))\ar[d]_{pr^*}
\ar[r]^(.55){\epsilon^*}&\Hom_{kF(\C)}(\P_*,\k)\ar[d]^{t^*}\\
\Hom_{k\C}(LK_{pr}LK_{\tau}(\P_*),LK_{pr}LK_{\tau}(\k))\ar[r]_(.55){\cong}&\Hom_{k\C}(LK_t(\P_*),LK_t(\k)).}
$$
Finally we show $\epsilon^*$ is a ring homomorphism. Since $\k =
\Res_{\tau}\k$, we get
$$
{\Ext}^*_{kF(\C)}(\k,\k)\cong{\Ext}^*_{kF(\C)}(\k,{\Res}_{\tau}\k)\cong{\Ext}^*_{k\C^e}({LK}_{\tau}(\k),\k)\cong{\Ext}^*_{k\C^e}(k\C,\k).
$$
It implies the cup product in the Hochschild cohomology ring
$$
{\Ext}^*_{k\C^e}(k\C,k\C){\otimes}_k{\Ext}^*_{k\C^e}(k\C,k\C){\buildrel{\smile}\over{\to}}
{\Ext}^*_{k\C^e}(k\C,k\C)={\Ext}^*_{k\C^e}(k\C,k\C{\otimes}_{k\C}k\C)
$$
is compatible with the cup product in the ordinary cohomology ring
since we have the following commutative diagram
$$
\xymatrix{\Ext^*_{k\C^e}(k\C,k\C)\otimes_k\Ext^*_{k\C^e}(k\C,k\C)\ar[r]^(.55){\smile}\ar@{>>}[d]_{\epsilon^*\otimes_k
\epsilon^*}
&\Ext^*_{k\C^e}(k\C,k\C\otimes_{k\C}k\C)\ar@{=}[r]& \Ext^*_{k\C^e}(k\C,k\C)\ar@{>>}[d]^{\epsilon^*}\\
\Ext^*_{k\C^e}(k\C,\k)\otimes_k\Ext^*_{k\C^e}(k\C,\k)\ar[r]^(.55){\smile}
&
\Ext^*_{k\C^e}(k\C,\k\otimes_{k\C}\k)\ar@{=}[r]&\Ext^*_{k\C^e}(k\C,\k).}
$$
Thus
$$
\epsilon^* :
{\Ext}^*_{k\C^e}(k\C,k\C)\twoheadrightarrow{\Ext}^*_{kF(\C)}(\k,\k)
$$
is a left inverse of $\tau^*$.
\end{proof}
\end{lem}

From the proof of last lemma, we have
$$
{\Ext}^*_{k\C^e}(k\C,M)\cong{\Ext}^*_{kF(\C)}(\k,{\Res}_{\tau}M)
$$
for any functor $M \in k\C^e$-mod. This is not necessarily true for any $M\in k\C^e$-mod
as $\tau : F(\C) \to \C^e$ does not always induce an algebra homomorphism hence the restriction
on $M$ may not make sense. Together with our earlier discussion, we
have the following formula for computation. Since we showed
$\Res_{\tau}(k\C)=T \oplus N_{\C}$ with $T\cong \k$, we may use the
decomposition to compute the Hochschild cohomology ring when the
structure of $N_{\C}$ is understood.

\begin{prop} Let $\C$ be a small category and $k$ a field. For any functor $M \in k\C^e$-mod, we have
$$
{\Ext}^*_{k\C^e}(k\C,M)\cong{\Ext}^*_{kF(\C)}(\k,{\Res}_{\tau}M).
$$
In particular we have
$$
{\Ext}^*_{k\C^e}(k\C,\k)\cong{\Ext}^*_{kF(\C)}(\k,\k)\cong{\Ext}^*_{k\C}(\k,\k),
$$
and
$$
{\Ext}^*_{k\C^e}(k\C,k\C)\cong{\Ext}^*_{kF(\C)}(\k,\k)\oplus{\Ext}^*_{kF(\C)}(\k,N_{\C})\cong{\Ext}^*_{k\C}(\k,\k)\oplus{\Ext}^*_{kF(\C)}(\k,N_{\C}),
$$
where $N_{\C}$ is the submodule of $\Res_{\tau}(k\C)\in kF(\C)$-mod which as a
functor takes the following value
$$
N_{\C}([\alpha])=k\{\beta - \gamma \bigm{|} \beta, \gamma \in
{\Hom}_{\C}(y,x) \},
$$
if $[\alpha]\in\Ob F(\C)$ and $\alpha\in\Hom_{\C}(y,x)$.
\end{prop}

Note that when $\C$ is a finite abelian group, we obtain Holm's
isomorphism \cite{Ho, CS}
$$
{\Ext}^*_{k\C^e}(k\C,k\C)\cong{\Ext}^*_{kF(\C)}(\k,{\Res}_{\tau}(k\C))\cong
k\C\otimes_k {\Ext}^*_{k\C}(\k,\k).
$$
In Section 3 we will compute some further examples of Hochschild
cohomology rings, using the above formula.

\subsection{EI-categories} A small category is EI if every endomorphism
is an isomorphism, and is finite if the morphism set is finite. The reader
is referred to \cite{We, X1} for a general description of the representation
and ordinary cohomology theory of finite EI-categories. In this subsection
we always assume $\C$ is a finite EI-category. The finiteness condition implies all
$k\C$-modules are functors, while the EI-condition implies that $x\cong x'$ in $\Ob\C$
if both $\Hom_{\C}(x,x')$ and $\Hom_{\C}(x',x)$ are non-empty. The EI-condition allows
us to give a partial order on the set of isomorphism classes of objects in $\Ob\C$
and hence a natural filtration to each functor in $k\C$-mod with respect to the partial order.
The simple and (finitely generated) projective $k\C$-modules have been classified by L\"uck \cite{Lu}.

For future reference, we quote the following
result \cite{X1}: let $\C$ be a finite EI-category and $M, N \in k\C$-mod.
An object $x\in\Ob\C$ is called $M$-minimal if $M(x)\ne 0$ and there is no
object $y \in \Ob\C$ such that $\Hom_{\C}(y,x)\ne\emptyset$ and $M(y)\ne 0$.
If the $M$-minimal objects are $x_1,\cdots,x_n\in\Ob\C$, and
$X_M$ is the full subcategory of $\C$ consisting of all $M$-minimal objects, then
$$
{\Ext}^*_{k\C}(M,N)\cong{\Ext}^*_{kX_M}(M,N),
$$
given that $N$ as a functor takes non-zero values only at objects
in $X_M$. This isomorphism will be used in this subsection as well as in the next section where
we compute some Hochschild cohomology rings.

Suppose $\A$ is the full subcategory of $\C$ which consists of all objects and all
isomorphisms in $\C$. The category $\A$ is a disjoint
union of finitely many finite groups. Its category algebra
$k\A=\oplus_{x\in\Ob\C}k\Aut_{\C}(x)$ is a $k\C^e$-module, and is a
quotient of $k\C$, with kernel written as $ker$. Considered as a
functor $ker \subset k\C$ takes non-zero values at $(x,y)$ for which
there exists a $\C$-morphism from $y$ to $x$ and $x\not\cong y$.

The short exact sequence of $k\C^e$-modules
$$
0 \to ker \to k\C {\buildrel{\pi}\over{\to}} k\A \to 0
$$
induces a long exact sequence
$$
{\cdots\to{\Ext}^n_{k\C^e}(k\C,ker)\to{\Ext}^n_{k\C^e}(k\C,k\C)
{\buildrel{\tilde{\pi}}\over{\to}}{\Ext}^n_{k\C^e}(k\C,k\A){\buildrel{\eta}\over{\to}}{\Ext}^{n+1}_{k\C^e}(k\C,ker)\to\cdots.}
$$
By the previously quoted result from \cite{X1}, one can see ${\Ext}^*_{k\C^e}(k\C,k\A)$ is naturally
isomorphic to
$$
{\Ext}^*_{k\A^e}(k\A,k\A),
$$
which is isomorphic to the direct sum of the Hochschild cohomology
rings of the automorphism groups of objects in $\C$:
$\oplus_{x\in\Ob\C}\Ext^*_{k\Aut_{\C}(x)^e}(k\Aut_{\C}(x),k\Aut_{\C}(x))$.
The following map will still be written as $\tilde{\pi}$
$$
\tilde{\pi} : {\Ext}^*_{k\C^e}(k\C,k\C)\to{\Ext}^*_{k\A^e}(k\A,
k\A).
$$
We show $\tilde{\pi}$ can be identified with the algebra
homomorphism induced by $-\otimes_{k\C}k\A$
$$
\phi_{k\A} : {\Ext}^*_{k\C^e}(k\C, k\C) \to {\Ext}^*_{k\C^e}(k\A,
k\A)\cong{\Ext}^*_{k\A^e}(k\A,k\A).
$$
Hence we do not need to distinguish the maps $\phi_{k\A}$ and
$\tilde{\pi}$.

\begin{lem} The following diagram is commutative
$$
\xymatrix{{\Ext}^*_{k\C^e}(k\C,k\C) \ar[r]^{\tilde{\pi}} \ar[d]_{\phi_{k\A}} & {\Ext}^*_{k\C^e}(k\C,k\A) \ar[d]^{\cong}\\
{\Ext}^*_{k\C^e}(k\A,k\A) \ar[r]_{\cong}
&{\Ext}^*_{k\A^e}(k\A,k\A).}
$$

\begin{proof}[Proof.] This can be seen on the cochain level. Suppose $\R_* \to k\C \to 0$
is the minimal projective resolution of the $k\C^e$-module $k\C$.
Then $\Ext^*_{k\C^e}(k\C,k\C)$ is the homology of the cochain
complex $\Hom_{k\C^e}(\R_*,k\C)$. The tensor product
$-\otimes_{k\C}k\A$ induces a map
$$
{\Hom}_{k\C^e}(\R_*,k\C)\to{\Hom}_{k\C^e}(\R_*\otimes_{k\C}k\A,k\C\otimes_{k\C}k\A)
\cong{\Hom}_{k\C^e}(\R_*\otimes_{k\C}k\A,k\A),
$$
which gives rise to $\phi_{k\A}$. On the other hand $\tilde{\pi}$ is
given by
$$
{\Hom}_{k\C^e}(\R_*,k\C)\to{\Hom}_{k\C^e}(\R_*,k\A)\cong{\Hom}_{k\A^e}({\Res}_{\C,\A}(\R_*),k\A),
$$
where ${\Res}_{\C,\A}(\R_*)$ is the restriction of $\R_*$ along the inclusion
$\A\hookrightarrow\C$ and is the minimal projective resolution of
the $k\A^e$-module $k\A$. But
$$
{\Hom}_{k\C^e}(\R_*\otimes_{k\C}k\A,k\A)\cong{\Hom}_{k\A^e}(\R_*\otimes_{k\C}k\A,k\A)
\cong{\Hom}_{k\A^e}({\Res}_{\C,\A}(\R_*),k\A).
$$
\end{proof}
\end{lem}

We have the following commutative diagram, involving four cohomology
rings.

\begin{thm} Let $\C$ be a finite EI-category and $k$ a field. Then we have the following
commutative diagram
$$
\xymatrix{{\Ext}^*_{k\C^e}(k\C,k\C) \ar[rr]^{\phi_{k\A}=\tilde{\pi}} \ar[d]_{\phi_{\C}} && {\Ext}^*_{k\A^e}(k\A,k\A) \ar[d]^{\phi_{\A}}\\
{\Ext}^*_{k\C}(\k,\k) \ar[rr]_{\Res_{\C,\A}} &&
{\Ext}^*_{k\A}(\k,\k).}
$$

\begin{proof}[Proof.] As usual, we prove it on the cochain level. Let $\R_* \to k\C \to 0$ be the
minimal projective resolution of the $k\C^e$-module $k\C$. Then we have the following
commutative diagram
$$
\xymatrix{\Hom_{k\C^e}(\R_*,k\C) \ar[r] \ar[d]& \Hom_{k\C^e}(\R_*\otimes_{k\C}k\A,k\C\otimes_{k\C}k\A) \ar[d]\\
\Hom_{k\C}(\R_*\otimes_{k\C}\k,k\C\otimes_{k\C}\k) \ar[r] \ar[d] & \Hom_{k\A}(\R_*\otimes_{k\C}k\A\otimes_{k\A}\k,k\C\otimes_{k\C}k\A\otimes_{k\A}\k)\ar[d]\\
\Hom_{k\C}(\R'_*,\k) \ar[r] & \Hom_{k\A}(\R'_*,\k) \ar[d] \\
& \Hom_{k\A}(\R''_*,\k),}
$$
in which $\R'_* \to \k \to 0$ and $\R''_* \to \k \to 0$ are the
projective resolutions of $k\C$- and $k\A$-modules
satisfying the following commutative diagrams of $k\C$-modules and
$k\A$-modules, respectively,
$$
\xymatrix{\R'_*\ar[r]\ar[d]& \k \ar[r]\ar[d]^{\cong}& 0&&\R''_*\ar[r]\ar[d]& \k \ar[r]\ar[d]^{=}& 0\\
\R_*\otimes_{k\C}\k \ar[r] & k\C\otimes_{k\C}\k \ar[r] &
0&\mbox{and}&\R'_* \ar[r] & \k \ar[r] & 0.}
$$
In the main diagram, upper left cochain complex computes
$\Ext^*_{k\C^e}(k\C,k\C)$, upper right corner computes
$\Ext^*_{k\A^e}(k\A,k\A)$, lower left corner computes
$\Ext^*_{k\C}(\k,\k)$ and lower right corner computes
$\Ext^*_{k\A}(\k,\k)$. Hence our statement follows.
\end{proof}
\end{thm}

We note that in the theorem the category $\A$ may be replaced by any full subcategory of it.
Especially, we have a commutative diagram for each $\Aut_{\C}(x)\subset\A$
$$
\xymatrix{{\Ext}^*_{k\C^e}(k\C,k\C)
\ar[rr]^(.33){\phi_{k\Aut_{\C}(x)}} \ar[d]_{\phi_{\C}} &&
{\Ext}^*_{k\Aut_{\C}(x)^e}(k\Aut_{\C}(x),k\Aut_{\C}(x)) \ar[d]^{\phi_{\Aut_{\C}(x)}}\\
{\Ext}^*_{k\C}(\k,\k) \ar[rr]_{\Res_{\C,\Aut_{\C}(x)}} &&
{\Ext}^*_{k\Aut_{\C}(x)}(\k,\k).}
$$

\section{Examples of the Hochschild cohomology rings of categories}

In this section we calculate the Hochschild cohomology rings for
four finite EI-categories, with base field $k$ of characteristic 2.
In particular the first category gives rise to a counterexample
against the finite generation conjecture of the Hochschild
cohomology rings in \cite{SS}.

Since all of our four categories are finite EI-categories, for the reader's convenience
we give a description of the simple $k\C$-modules for a finite EI-category $\C$.
By \cite{Lu}, any simple $k\C$-module $S_{x,V}$ is indexed by the isomorphism class of
an object $x\in\Ob\C$ and a simple module $V$ of the automorphism group $\Aut_{\C}(x)$ of $x$.
As a functor, $S_{x,V}(y)\cong V$ if $y\cong x$ in $\Ob\C$ and $S_{x,V}(y)=0$ otherwise.

\subsection{The category $\E_0$}

In \cite{X2} we presented an example, by Aur\'elien Djament, Laurent
Piriou and the author, of the mod-2 ordinary cohomology ring of the
following category $\E_0$
$$
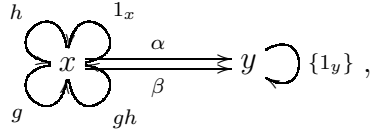
\xymatrix{x \ar@(r,u)[]_{1_x} \ar@(l,d)[]_g \ar@(u,l)[]_h
\ar@(d,r)[]_{gh} \ar@<2pt>[rr]^{\alpha} \ar@<-2pt>[rr]_{\beta}&& y
\ar@(ur,dr)[]^{\{1_y\}}} ,
$$
where $g^2=h^2=1_x, gh=hg, \alpha h = \beta g = \alpha$, and $\alpha
g = \beta h = \beta$. The ordinary cohomology ring $\Ext^*_{k\E_0}(\k,\k)$ is a subring of the polynomial ring
$\H^*(\Z_2\times\Z_2,k)\cong k[u,v]$, removing all $u^n, n\ge 1$, and
their scalar multiples.
It has no nilpotents and is not finitely generated. By Theorem
2.3.4, it implies that the Hochschild cohomology ring
$\Ext^*_{k\E_0^e}(k\E_0,k\E_0)$ is not finitely generated either,
which gives a counterexample against the conjecture in \cite{SS}. We
compute its Hochschild cohomology ring using Proposition 2.3.5.

The category of factorizations in $\E_0$, $F(\E_0)$, has the
following shape
$$
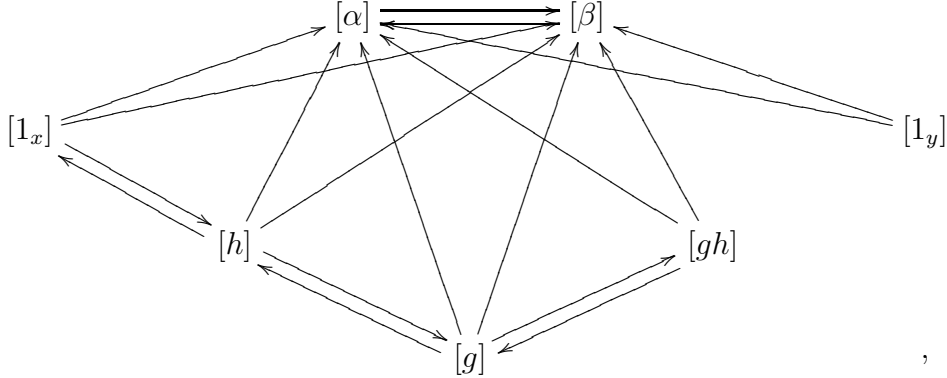
\xymatrix{&&&[\alpha] \ar@<.5ex>[rr]& & [\beta] \ar@<.5ex>[ll]&&&\\
[1_x] \ar@<.5ex>[drr] \ar[urrr] \ar[urrrrr]&&&&&&&& [1_y] \ar[ulll] \ar[ulllll]\\
&& [h] \ar@<.5ex>[ull] \ar[uurrr] \ar[uur] \ar@<.5ex>[drr]&&&& [gh] \ar[uulll] \ar[uul] \ar@<.5ex>[dll]&&\\
&&&& [g] \ar[uuul] \ar[uuur] \ar@<.5ex>[ull] \ar@<.5ex>[urr]&&&&,}
$$
in which $[1_x]\cong [h]\cong [g] \cong [gh]$ and
$[\alpha]\cong[\beta]$. For the purpose of computation, we use the
skeleton $F'(\E_0)$ of $F(\E_0)$ (which is equivalent to $F(\E_0)$ hence the two
category algebras and their module categories are Morita equivalent)
$$
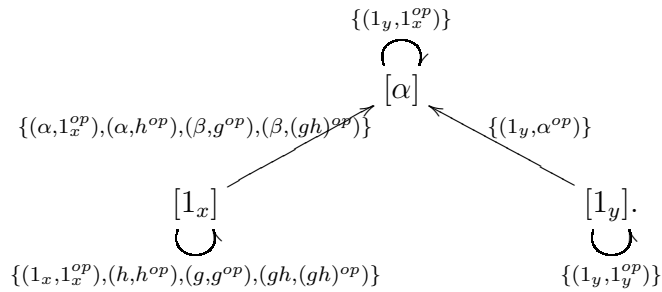
\xymatrix{&&[\alpha] \ar@(ul,ur)^{\{(1_y,1_x^{op})\}}&& \\
[1_x] \ar[urr]^{\{(\alpha, 1_x^{op}), (\alpha,h^{op}), (\beta,g^{op}),
(\beta, (gh)^{op})\}} \ar@(dl,dr)_{\{(1_x,1_x^{op}), (h,h^{op}),
(g,g^{op}), (gh,(gh)^{op})\}} &&&& [1_y]. \ar[ull]_{\{(1_y,\alpha^{op})\}}
\ar@(dl,dr)_{\{(1_y,1_y^{op})\}}}
$$
In the above category, next to each arrow is the set of
homomorphisms in $F'(\E_0)$ from one object to another. The module
$N_{\E_0} \in kF'(\E_0)$-mod (see Proposition 2.3.5) takes the
following values
$$
\begin{array}{lllllll}
N_{\C}([1_x])&=& k\{1_x + h, g + gh, 1_x + g\} &,& N_{\C}([h])&=&k\{1_x + h, g + gh, 1_x + g\},\\
N_{\C}([g])&=&k\{1_x + h, g + gh, 1_x + g\}&,& N_{\C}([gh])&=& k\{1_x + h, g + gh, 1_x + g\},\\
N_{\C}([\alpha])&=& k\{\alpha + \beta\}&,& N_{\C}([\beta])&=&k\{\alpha + \beta\},\\
N_{\C}([1_y])&=& 0. &&&&
\end{array}
$$
Thus $N_{\E_0}=S_{[1_x],k(1_x + h)}\oplus S_{[1_x], k(g + gh)}
\oplus \k'_{1_x+g}$, where $S_{[1_x], k(1_x+h)}$ and
$S_{[1_x],k(g+gh)}$ are simple $kF'(\E_0)$-modules such that
$S_{[1_x],k(1_x + h)}([1_x])=k(1_x + h)$ and $S_{[1_x],
k(g+gh)}([1_x])=k(g+gh)$, and $\k'_{1_x+g}$ is a $kF'(\E_0)$-module
such that $\k'_{1_x+g}([1_x])=k(1_x+g)$,
$\k'_{1_x+g}([\alpha])=k(\alpha+\beta)$ and $\k'_{1_x+g}([1_y])=0$.
Note that $S_{[1_x],k(1_x + h)}([1_x])=k(1_x+h)$, $S_{[1_x],
k(g+gh)}([1_x])=k(g+gh)$ and $\k'_{1_x+g}([1_x])=k(1_x+g)$ are all
isomorphic to the trivial $k\Aut_{F'(\E_0)}([1_x])$-module, and have
the same trivial ring structure in the sense that the product of any
two elements is zero. Hence we have (along with the result quoted in
Section 2.4, paragraph two)
$$
{\Ext}^*_{kF'(\E_0)}(\k,S_{[1_x],k(1_x+h)})\cong
k(1_x + h)\otimes_k{\Ext}^*_{k\Aut_{F'(\E_0)}([1_x])}(k,k)
$$
and
$$
{\Ext}^*_{kF'(\E_0)}(\k,S_{[1_x],k(g+gh)})\cong k(g +
gh)\otimes_k{\Ext}^*_{k\Aut_{F'(\E_0)}([1_x])}(k,k)
$$
as rings, in which $k(1_x+h)$ and $k(g+gh)$ are concentrated in
degree zero in each ring. From the structure of $F(\E_0)$, one has
$\Aut_{F'(\E_0)}([1_x])\cong\Z_2\times\Z_2$.

For computing $\Ext^*_{kF'(\E_0)}(\k,\k'_{1_x + g})$, we use the
following short exact sequence of $kF(\E_0)$-modules
$$
0 \to \k'_{1_x+g} \to \k \to S_{[1_y],k} \to 0.
$$
It induces a long exact sequence in which one can find
$\Ext^0_{kF'(\E_0)}(\k,S_{[1_y],k})=k$ and
$\Ext^n_{kF'(\E_0)}(\k,S_{[1_y],k})=0$ if $n \ge 1$. Thus
$\Ext^0_{kF'(\E_0)}(\k,\k'_{1+g})=0$ while
$\Ext^n_{kF'(\E_0)}(\k,\k)\cong\Ext^n_{kF'(\E_0)}(\k,\k'_{1+g})$ for
each $n\ge 1$. Hence as a ring
$$
{\Ext}^*_{kF'(\E_0)}(\k,\k'_{1_x+g})\cong
k(1_x+g)\otimes_k{\Ext}^{*>0}_{kF'(\E_0)}(\k,\k)\cong
k(1_x+g)\otimes_k{\Ext}^{*>0}_{k\E_0}(\k,\k).
$$
All in all, we have
$$
{\Ext}^0_{k\E_0^e}(k\E_0,k\E_0)\cong{\Ext}^0_{k\E_0}(\k,\k)\oplus
k(1_x + h)\oplus k(g + gh),
$$
and if $n\ge 1$
$$
\begin{array}{ll}
&{\Ext}^n_{k\E_0^e}(k\E_0,k\E_0)\\
\cong & {\Ext}^n_{k\E_0}(\k,\k) \oplus
\{k(1_x+g)\otimes_k\Ext^n_{k\E_0}(\k,\k)\}\\
& \oplus \{k(1_x + h)\otimes_k\Ext^*_{k(\Z_2\times\Z_2)}(k,k)\}
\oplus \{k(g + gh)\otimes_k\Ext^n_{k(\Z_2\times\Z_2)}(k,k)\}.
\end{array}
$$
Combining all the information we obtained, the surjective ring
homomorphism
$$
\phi_{\E_0} : {\Ext}^*_{k\E_0^e}(k\E_0,k\E_0) \twoheadrightarrow
{\Ext}^*_{k\E_0}(\k,\k)
$$
has its kernel consisting of all nilpotents. Consequently this
Hochschild cohomology ring modulo nilpotents is not finitely
generated, against the finite generation conjecture in \cite{SS}. We
comment that the category algebra $k\E_0$ is not a
self-injective algebra (hence is not Hopf, by \cite{LS}). Nicole Snashall points out to the author
that this algebra is Koszul since both $k\E_0$ and
$\Ext^*_{k\E_0}(\overline{k\E_0},\overline{k\E_0})$ as graded
algebras are generated in degrees zero and one, where
$\overline{k\E_0}=k\E_0/\Rad(k\E_0)\cong S_{x,k}\oplus S_{y,k}$.

\subsection{The category $\E_1$} The following category $\E_1$ has a terminal
object and hence is contractible:
$$
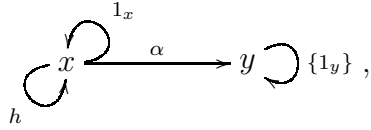
\xymatrix{x \ar@(r,u)[]_{1_x} \ar@(l,d)[]_h \ar[rr]^{\alpha} && y
\ar@(ur,dr)[]^{\{1_y\}}} ,
$$
where $h^2=1_x$ and $\alpha h=\alpha$. The contractibility implies
the ordinary cohomology ring is simply the base field $k$. In this
case $F(\E_1)$ is the following category
$$
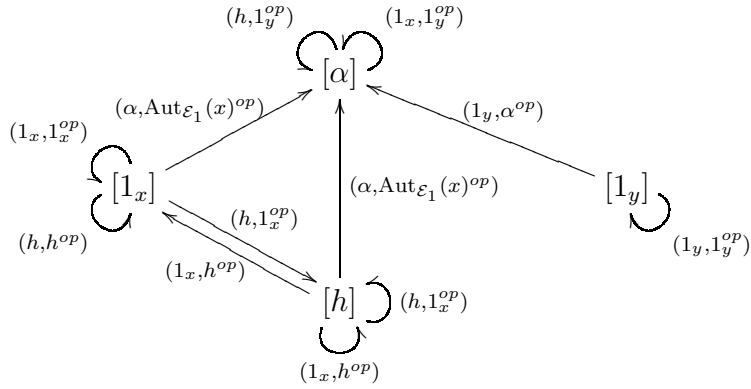
\xymatrix{&& [\alpha] \ar@(r,u)[]_{(1_x,1_y^{op})} \ar@(u,l)[]_{(h,1_y^{op})}&&&\\
[1_x] \ar@(u,l)[]_{(1_x,1_x^{op})} \ar@(l,d)[]_{(h,h^{op})}
\ar[urr]^{(\alpha,\Aut_{\E_1}(x)^{op})}
\ar@<.5ex>[drr]^{(h,1_x^{op})}&& &&& [1_y]
\ar[ulll]_{(1_y,\alpha^{op})} \ar@(r,d)[]^{(1_y,1_y^{op})}\\
&& [h] \ar@(dl,dr)[]_{(1_x,h^{op})} \ar@(dr,ur)[]_{(h,1_x^{op})}
\ar[uu]_{(\alpha,\Aut_{\E_1}(x)^{op})}
\ar@<.5ex>[ull]^{(1_x,h^{op})}&&&}
$$

We calculate its Hochschild cohomology ring. By proposition 2.3.5,
we only need to compute $\Ext^*_{kF(\E_1)}(\k,N_{\E_1})$, where
$N_{\E_1}$ has the following value at objects of $F(\E_1)$
$$
\begin{array}{lllllll}
N_{\E_1}([1_x])&=& k\{1_x + h\} &,& N_{\E_1}([h])&=&k\{1_x + h\},\\
N_{\E_1}([1_y])&=& 0 &,& N_{\E_1}([\alpha])&=&0.
\end{array}
$$
One can easily see that $N_{\E_1}=S_{[1_x], k(1_x+h)}$ is a simple
module of dimension one with a specified value $k(1_x+h)$ at
$[1_x]$. Since $[1_x] \cong [h] \in \Ob F(\E_1)$ are minimal
objects, using quoted result in Section 2.4 paragraph two, we get
$$
{\Ext}^*_{kF(\E_1)}(\k,N_{\E_1})\cong{\Ext}^*_{k\Aut_{F(\E_1)}([1_x])}(k,k(1_x
+ h))\cong k(1_x + h)\otimes_k{\Ext}^*_{k\Z_2}(k,k),
$$
which is isomorphic to $k(1_x + h)\otimes_k k[u]$. Here $k[u]$ is a
polynomial algebra with an indeterminant $u$ at degree one and
$k(1_x + h)$ is at degree zero. Thus
$$
{\Ext}^*_{k\E_1^e}(k\E_1,k\E_1)\cong{\Ext}^*_{k\E_1}(\k,\k)\oplus{\Ext}^*_{kF(\E_1)}(\k,N_{\E_1})\cong
k\oplus \{k(1_x + h)\otimes_k k[u]\}.
$$
The kernel of $\phi_{\E_1}$ consists of all nilpotents in the
Hochschild cohomology ring.

\subsection{The category $\E_2$} The following category has its classifying
space homotopy equivalent to the join,
$B\Z_2*B\Z_2=\Sigma(B\Z_2\wedge B\Z_2)=\Sigma[B(\Z_2\times
\Z_2)/(B\Z_2\vee B\Z_2)]$, of the classifying spaces of the two
automorphism groups:
$$
\xymatrix{x \ar@(r,u)[]_{1_x} \ar@(l,d)[]_h \ar[rr]^{\alpha} && y
\ar@(r,u)[]_{1_y} \ar@(l,d)[]_g} ,
$$
where $h^2=1_x$, $\alpha h = \alpha =g\alpha$ and $g^2=1_y$. As
direct consequences, its ordinary cohomology groups are equal to $k,
0, 0$ at degrees zero, one and two, and $k^{n-2}$ at each degree
$n\ge 3$, and furthermore the cup product in this ring is trivial
\cite{X2}. We compute its Hochschild cohomology ring. The category
$F(\E_2)$ is as follows
$$
\xymatrix{&&&& [\alpha] \ar@(ul,ur)^{(\Aut_{\E_2}(x),\Aut_{\E_2}(y)^{op})}&&&&\\
[1_x]
\ar[urrrr]^{(\alpha,\Aut_{\E_2}(x)^{op})}\ar@(u,l)[]_{(1_x,1_x^{op})}
\ar@(l,d)[]_{(h,h^{op})} \ar@<.5ex>[drr]
&&&&&&&& [g] \ar[ullll]_{(\Aut_{\E_2}(y),\alpha^{op})} \ar@(r,u)[]_{(1_y,g^{op})} \ar@(r,d)[]^{(g,1_y^{op})} \ar@<.5ex>[dll]\\
&& [h] \ar[uurr]^{(\alpha,\Aut_{\E_2}(x)^{op})}
\ar@(l,d)[]_{(1_x,h^{op})} \ar@(d,r)[]_{(h,1_x^{op})}
\ar@<.5ex>[ull]&&&& [1_y] \ar[uull]_{(\Aut_{\E_2}(y),\alpha^{op})}
\ar@(l,d)[]_{(1_y,1_y^{op})} \ar@(r,d)[]^{(g,g^{op})}
\ar@<.5ex>[urr] &&}
$$

By Proposition 2.3.5, we need to compute
$\Ext^*_{k\E_2}(\k,N_{\E_2})$. In this case we have
$$
\begin{array}{lllllll}
N_{\E_2}([1_x])&=& k\{1_x + h\} &,& N_{\E_2}([h])&=&k\{1_x + h\},\\
N_{\E_2}([1_y])&=& k\{1_y + g\} &,& N_{\E_2}([g])&=& k\{1_y + g\},\\
N_{\E_2}([\alpha])&=&0. &&&&
\end{array}
$$
It means $N_{\E_2}=S_{[1_x], k(1_x+h)}\oplus S_{[1_y],k(1_y+g)}$ and
thus by Proposition 2.2.5
$$
\begin{array}{ll}
{\Ext}^*_{k\E_2}(\k,N_{\E_2})&\cong{\Ext}^*_{k\Aut_{F(\E_2)}([1_x])}(k,k(1_x
+ h))\oplus{\Ext}^*_{k\Aut_{F(\E_2)}([1_y])}(k,k(1_y + g))\\
&\cong \{k(1_x +h)\otimes_k{\Ext}^*_{k\Z_2}(k,k)\}\oplus
\{k(1_y+g)\otimes_k{\Ext}^*_{k\Z_2}(k,k)\}.
\end{array}
$$
Hence
$$
{\Ext}^*_{k\E_2^e}(k\E_2,k\E_2)\cong{\Ext}^*_{k\E_2}(\k,\k)\oplus\{k(1_x
+ h)\otimes_k k[u]\}\oplus\{k(1_y + g)\otimes_k k[v]\},
$$
where $k[u]$ and $k[v]$ are two polynomial algebras with
indeterminants in degree one. Both the Hochschild and ordinary
cohomology rings modulo nilpotents are isomorphic to the base field
$k$.

\subsection{The category $\E_3$} The following category has a
classifying space homotopy equivalent to that of
$\Aut_{\E_3}(x)\cong\Z_2$ (by Quillen's Theorem A \cite{Qu}, or
see \cite{X1})
$$
\xymatrix{x \ar@(r,u)[]_{1_x} \ar@(l,d)[]_h \ar@<2pt>[rr]^{\alpha}
\ar@<-2pt>[rr]_{\beta}&& y \ar@(ur,dr)[]^{\{1_y\}}} ,
$$
where $h^2=1_x$ and $\alpha h=\beta$. We compute its Hochschild
cohomology ring. The category $F(\E_3)$ is as follows (not all
morphisms are presented since only its skeleton is needed)
$$
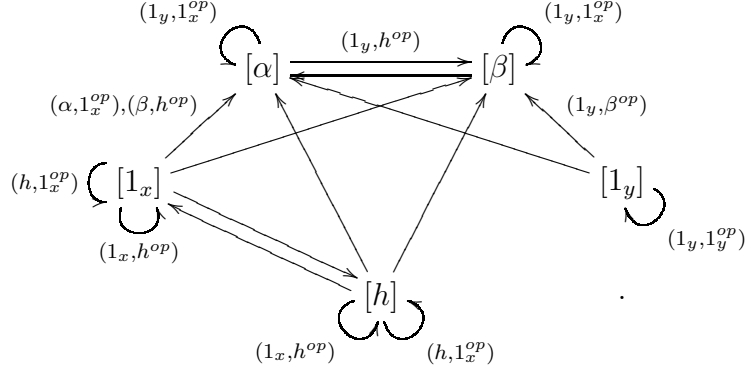
\xymatrix{& [\alpha] \ar@<.5ex>[rr]^{(1_y,h^{op})} \ar@(u,l)[]_{(1_y,1_x^{op})} & & [\beta] \ar@<.5ex>[ll] \ar@(u,r)[]^{(1_y,1_x^{op})}&\\
[1_x] \ar@<.5ex>[drr] \ar[ur]^{(\alpha,1_x^{op}),(\beta,h^{op})}
\ar[urrr] \ar@(dl,dr)[]_{(1_x,h^{op})}
\ar@(ul,dl)[]_{(h,1_x^{op})} && && [1_y] \ar[ulll] \ar[ul]_{(1_y,\beta^{op})} \ar@(r,d)^{(1_y,1_y^{op})}\\
&& [h] \ar@<.5ex>[ull] \ar[uul] \ar[uur] \ar@(l,d)[]_{(1_x,h^{op})}
\ar@(d,r)[]_{(h,1_x^{op})} &&.}
$$
The module $N_{\E_3}$ takes the following values
$$
\begin{array}{lllllll}
N_{\E_2}([1_x])&=& k\{1_x + h\} &,& N_{\E_2}([h])&=&k\{1_x + h\},\\
N_{\E_2}([1_y])&=& 0 &,& N_{\E_2}([\alpha])&=& k\{\alpha + \beta\},\\
N_{\E_2}([\alpha])&=& k\{\alpha + \beta\}. &&&&
\end{array}
$$
Thus $N_{\E_3}$ fits into the following short exact sequence of
$kF(\E_3)$-modules
$$
0\to N_{\E_3} \to \k \to S_{[1_y],k} \to 0.
$$
Just like in our first example, using the long exact sequence coming
from it, we know $\Ext^0_{k\E_2}(\k,N_{\E_3})=0$ and
$\Ext^{*>0}_{k\E_2}(\k,N_{\E_3})\cong k(1_x + h)\otimes_k
\Ext^{*>0}_{kF(\E_3)}(\k,\k)\cong k(1_x + h)\otimes_k
\Ext^{*>0}_{k\E_3}(\k,\k)$. Hence
$$
{\Ext}^*_{k\E_3^e}(k\E_3,k\E_3)\cong{\Ext}^*_{k\E_3}(\k,\k)\oplus
\{k(1_x + h)\otimes_k {\Ext}^{*>0}_{k\E_3}(\k,\k)\}.
$$
The kernel of $\phi_{\E_3}$ contains all nilpotents in the
Hochschild cohomology ring.


\begin{thebibliography}{50}
\bibitem{BLO} C. Broto, R. Levi and B. Oliver, The homotopy theory of
fusion systems, J. Amer. Math. Soc. 16 (2003), 779-856.
\bibitem{Be} D. Benson, Representations and Cohomology II, Cambridge
Studies in Adv. Math. 31, Cambridge University Press, 1998.
\bibitem{BW} H.-J. Baues and G. Wirsching, Cohomology of small
categories, J. Pure Appl. Algebra 38 (1985), 187-211.
\bibitem{CS} C. Cibils and A. Solotar, Hochschild cohomology algebra of abelian groups,
Arch. Math. 68 (1997), 17-21.
\bibitem{FS} E. Friedlander and A. Suslin, Cohomology of finite group
schemes over a field, Invent. Math. 127 (1997) no. 2, 209-270.
\bibitem{Ge} M. Gerstenhaber, The cohomology structure of an associative ring, Ann.
Math. 78 (1963), 267-288.
\bibitem{GS} M. Gerstenhaber and S. Shack, Simplicial cohomology is
Hochschild cohomology, J. Pure Appl. Algebra 30 (1983), 143-156.
\bibitem{GSS} E. Green, N. Snashall and $\O$. Solberg, The Hochschild cohomology ring
of a self-injective algebra of finite representation type, Proc. Amer. Math. Soc.
131 (2003), 3387-3393.
\bibitem{Ho} T. Holm, The Hochschild cohomology ring of a modular group
algebra: the commutative case, Comm. Algebra 24 (1996), 1957-1969.
\bibitem{HV} J. Hollender and R. Vogt, Modules of topological spaces, applications
to homotopy limits and $E_{\infty}$-structures, Arch. Math. 59 (1992) 115-129.
\bibitem{LS} R. Larson and M. Sweedler, An associative orthogonal bilinear form for Hopf algebras, Amer. J. Math. 91 (1969) 75-94.
\bibitem{L1} M. Linckelmann, Varieties in block theory, J. Algebra 215 (1999), 460-480.
\bibitem{L2} M. Linckelmann, Alperin's weight conjecture in terms
of Bredon equivariant cohomology, Math. Z. 250 (2005), 493-513.
\bibitem{Lo} J.-L. Loday, Cyclic Homology (Second Edition), Grundlehren der mathematischen Wissenschaften 301, Springer-Verlag Berlin, 1998.
\bibitem{Lu} W. L{\" u}ck, Transformation Groups and Algebraic K-Theory,
Lecture Notes in Math. 1408, Springer-Verlag Berlin, 1989.
\bibitem{ML} S. Mac Lane, Homology, Springer-Verlag Berlin, 1995.
\bibitem{Mi} B. Mitchell, Rings with several objects, Adv. in Math. 8 (1972), 1-161.
\bibitem{Qu} D. Quillen, Higher algebraic $K$-theory I, Lecture Notes
in Math. 341, 85-147, Springer-Verlag Berlin, 1973.
\bibitem{Sl} J. S{\l}omi\'nska, Homotopy colimits on EI-categories,
Lecture Notes in Mathematics 1408, Springer-Verlag Berlin (1989),
273-294.
\bibitem{SS} N. Snashall and $\O$. Solberg, Support varieties and
Hochschild cohomology rings, Proc. London Math. Soc., vol 88,
no. 3 (2004), 705-732.
\bibitem{Th} J. Th\'evenaz, G-algebras and Modular
Representation Theory, Oxford Mathematical Monographs, The Clarendon
Press, Oxford University Press, New York, 1995.
\bibitem{We} P. J. Webb, An introduction to the representations
and cohomology of categories, Group Representation Theory, 149-173,
EPFL Press, Lausanne, 2007.
\bibitem{X1} F. Xu, Representations of categories and their applications, J. Algebra 317 (2007), 153-183.
\bibitem{X2} F. Xu, On the cohomology rings of small categories, J. Pure Appl. Algebra, to appear.
\end{thebibliography}
\end{document}